\documentclass[12pt,a4paper]{amsart}
\usepackage{amsfonts}
\usepackage{amsthm}
\usepackage{amsmath}
\usepackage{amscd}
\usepackage[utf8]{inputenc}
\usepackage{t1enc}
\usepackage[mathscr]{eucal}
\usepackage{indentfirst}
\usepackage{graphicx}
\usepackage{graphics}
\usepackage{pict2e}
\usepackage{epic}
\usepackage{float}
\numberwithin{equation}{section}
\usepackage[margin=2.9cm]{geometry}
\usepackage{epstopdf} 
\usepackage{hyperref}
\hypersetup{
	colorlinks = true,
	urlcolor   = blue,
	citecolor  = blue,
} 
\usepackage{tikz-cd}
\usepackage{verbatim}
\usepackage{float}
\usepackage{makecell}
\begin{document}
	\newtheorem{de}{Definition}
	\newtheorem{ex}[de]{\emph{Example}}
	\newtheorem{thm}[de]{Theorem}
	\newtheorem{lemma}[de]{Lemma}
	\newtheorem{cor}[de]{Corollary}
	\newtheorem{con}[de]{Conjecture}
	\newtheorem{prop}[de]{Proposition}
	\title{An upper bound of the numbers of minimally intersecting filling coherent pairs}
	
	\author{Hong Chang}
	
	\address{
		\noindent
		Hong Chang, hchang24@buffalo.edu, Department of Mathematics, University at Buffalo--SUNY}

	
	\begin{abstract}
		Let $S_g$ denoting the genus $g$ closed orientable surface. An {\em origami} (or flat structure) on $S_g$ is obtained from a finite collection of unit Euclidean squares by gluing each right edge to a left one and each top edge to a bottom one.
		Coherent filling pairs of simple closed curves, $(\alpha,\beta)$ in $S_g$ are pairs for which their minimal intersection is equal to their algebraic intersection. And, a minimally intersecting filling of $(\alpha,\beta)$ in $S_g$ is a pair whose intersection number is the minimal among all filling pairs of $S_g$.
		A coherent pair of curves is naturally associated with an origami on $S_g$, and a minimally intersecting filling coherent pair of curves has the smallest number of squares in all origamis on $S_g$.
		Our main result introduces an algorithm to count the numbers of minimal filling pairs on $S_g$, and establish a new upper bound of this count using M\'enage Problem by \'Edouard Lucas in \cite{lucas}.
	\end{abstract}
	\maketitle
	\section{Introduction}
	Let $S_g$ denote the closed orientable surface of genus $g$. A finite collection of curves $\Gamma = {\gamma_1, ...\gamma_n}$ in
	pairwise minimal position on $S_g$ is said to fill if the complement $S_g\backslash\gamma$ is a disjoint union of topological
	disks. When each $\gamma_i$ is simple---no self intersections---and when n = 2, we call $\Gamma$ a filling pair.
	
	Let $\alpha$ and $\beta$ be a filling pair on $S_g$. The pair, $(\alpha, \beta)$, is a \emph{minimally intersecting} filling pair if the intersecting number $i(\alpha,\beta)$ is minimal among all filling pairs on $S_g$. 
	
	Aougab and Huang showed in \cite{AMN} that for all $g >2$ there exists filling pairs of curves whose intersection achieves the $2g-1$ minima, and when $g=2$, the minima is 4. In fact, the minima can be obtained such that the absolute value of the algebraic intersection of $\alpha$ and $\beta$ is equal to the geometric intersection---$\alpha$ and $\beta$ \emph{intersect coherently}.  So, \emph{minimally intersecting coherent filling  pairs} is a non-empty category for all $S_{g\ge 2}$.
	
	Let $u$ and $v$ be the isotopy classes of $\alpha$ and $\beta$, and let $\alpha$ and $\beta$ be the minimal intersection number over all isotopic
	representatives. We will abuse notation by having $i(u,v)=i(\alpha,\beta)$. For convenience, we will also say $(u,v)$ is a \emph{filling pair} when $(\alpha, \beta)$ is a filling pair. In this paper, when counting the number of filling pairs, we see $(u,v)$ as an ordered pair of classes.  Necessarily, curves and classes have no orientation.
	
	The \emph{mapping} class group of S, denoted $Mod(S)$, is the group of isotopy classes of orientation preserving self-homeomorphisms of S. The \emph{extended mapping class group} of S, denoted $Mod^{\pm}(S)$, is the group of isotopy classes of self-homeomorphisms of S, including the orientation-reversing ones. We say $\alpha$, $\beta$ and $\alpha'$, $\beta'$ are two filling pairs in the same $Mod(S)$ (or $Mod^{\pm}(S)$) orbit if there exist $g\in Mod(S)$ (or $Mod^{\pm}(S)$) such that $g(\alpha)=\alpha'$ and $g(\beta)=\beta'$.
	
	Before we go to the bounds, we give a definition for asymptotic notation ``$\sim''$:
	
	\begin{de}
	    If $f_n$ and $g_n$ are sequences of real valued function and $g_i\not=0$ for all $i$. We call $f_n\sim\ g_n$ if $\underset{n\rightarrow\infty}{\lim}f_n/g_n=1$.
	\end{de}

 Use $\mathcal{MI}$ for the set of all $Mod(S)$ orbit of minimally intersecting pairs.  Then the size of the set is $$f(g)\leq |\mathcal{MI}| \leq 2^{2g-2}(4g-5)(2g-3)!$$Where $f(g)\sim 3^{g/2}/g^2$. This result is provided in Aoubag--Huang \cite{AH}.

 Use $\mathcal{CMI}$ for $Mod^{\pm}(S)$ orbit of coherently minimally intersecting pairs. Then the size of the set is $$p(g)\leq |\mathcal{CMI}| \leq h(g)$$Where $p(g)=\begin{cases} (g-2)!  & \text{ if } g>2 \text{ odd} \\(g-5)(g-3)! & \text{ if } g>2 \text{ even} \end{cases}$ and $h(g)\sim\frac{(g-1)(2g-2)!}{e^2}$. The left side is provided in \cite{AMN} and the right side will be proved in this paper.
	\begin{thm} \label{main}
		Let $\mathcal{CMI}$ be $Mod^{\pm}(S)$ orbit of coherently minimally intersecting pairs on $S_g$, then $|\mathcal{CMI}| \leq h(g)$ where $$h(g)\sim\frac{(g-1)(2g-2)!}{e^2}$$.
	\end{thm}
	\subsection{Origamis}
	An \emph{origami}, or a square-tiled surface for a closed surface $S_{g \geq 2}$ is obtained from a finite number of Euclidean squares by gluing each right edge to a left one and each top edge to a bottom one. A \emph{$[1,1]$-origami} will be an origami that has exactly one horizontal annulus and one vertical annulus. In particular, we have the following well known result associating a $[1,1]$-origami to an origami pair of curves.
	
	\begin{thm}[Theorem 1.1 of \cite{CJM}]
	    A coherent filling pair of curves (origami pair of curves) naturally corresponds to an origami on $S_g$.
	\end{thm}
	
	As mentioned in \cite{AMN}, an abelian differential on a surface of genus $g > 1$ must have $2g - 2$ zeros, counted with multiplicity. In the moduli space of abelian differentials we denote the stratum $H(m_1, ..., m_n)$ of the moduli space consists of those abelian differentials which have $n$ zeros of degrees $m_1+...+m_n=2g-2$.
	
	Eskin-Okounkov \cite{Esk} and Zorich \cite{Zor} are using square-tiled surfaces (origamis) to calculate the volume of a stratum of abelian differentials. And for this reason, counting the number of origamis has connections to Teichm\"uller dynamics and to the study of flat surfaces. In this fashion of origamis, Theorem \ref{main} can be also expressed as:
	
	\emph{In the minimal stratum $H(2g-2)$, there exist at most $h(g)$ $[1,1]$-origamis.}
	
	\subsection{Construction of surface from polygons}
	Since $\alpha$ and $\beta$ is minimally filling, when $g\ge 3$, $S_g\backslash\alpha\cup\beta$ will consist of only one component, which is a $4(2g-1)$-gon with 2 copies of $\alpha$ and 2 copies of $\beta$ as the sides. 
	\begin{figure}[h]
		\scalebox{.5}{\includegraphics[angle=0,origin=c]{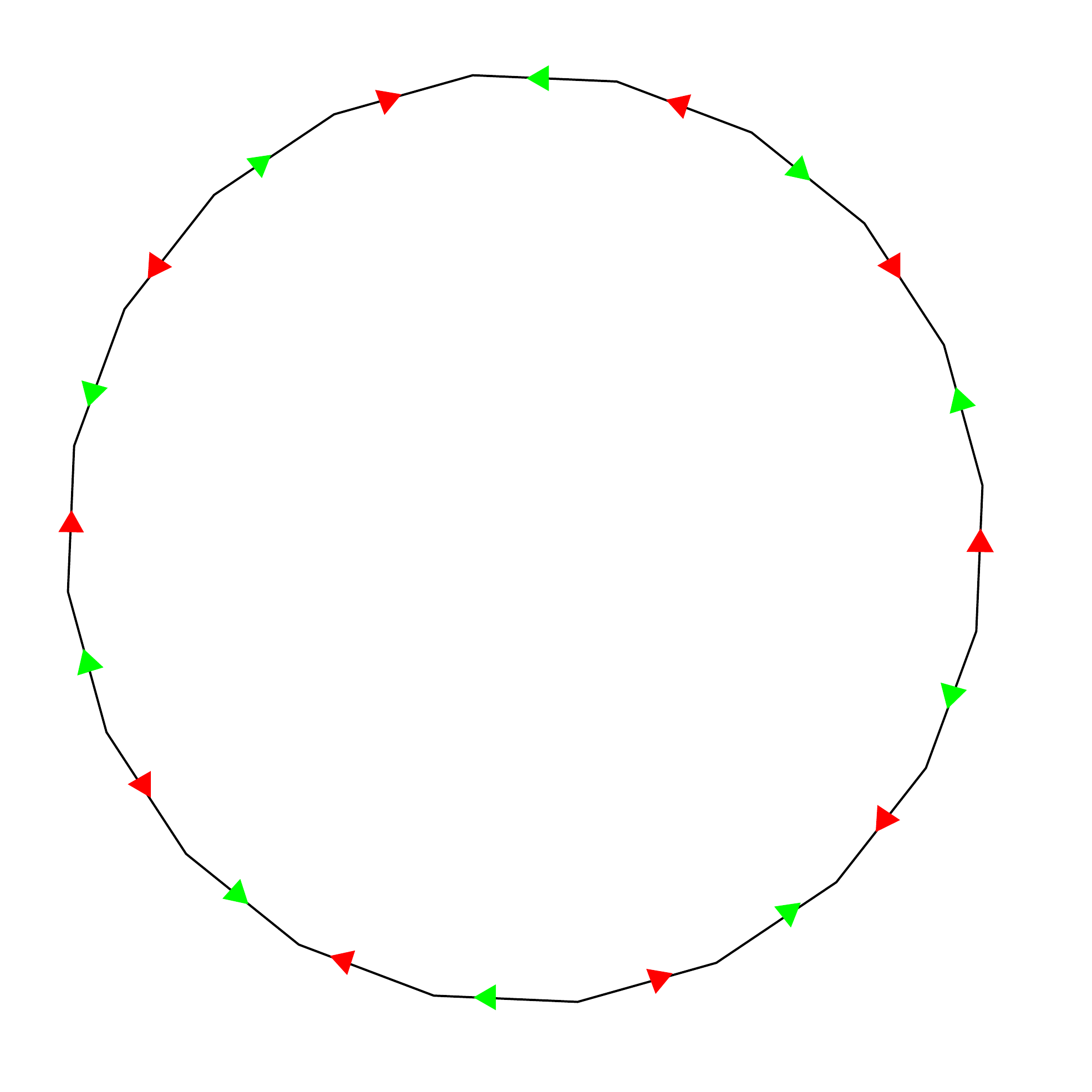}}
		\vspace{0cm}
		\caption{Coherent intersecting sides of $4(2g-1)$-gons, where $g$ is odd(left) and even(right). The red is $\alpha$ while green is $\beta$.}
		\label{gon}
	\end{figure}
	
	Please notice the edges for the $4(2g-1)$-gon are one-to-one correspond with subarcs of $\alpha$ and $\beta$, thus the horizontal and vertical lines of the origami.
	\subsection{M\'enage Problem}
	M\'enage Problem is asked by \`Edouard Lucas in \cite{lucas} as ``the number of different ways in which it is possible to seat a set of male-female couples at a round dining table so that men and women alternate and nobody sits next to his or her partner'', or by Peter Tait in \cite{tait} as ``numbers of arrangements such that there are of n letters, when A cannot be in the first or second place, B not in the second or third'' (such arrangements can be regarded as permutations called \emph{m\'enage permutation}). The first one will be equivalent to the second one if we already seat, say wives, on the table, which there are $2n!$ ways to do that.
	
	Gilbert, E. N. modified Tait's problem further in \cite{gen}. Let $P$ be a m\'enage permutation and $C$ a cyclic permutation $i\rightarrow i+a$ (mod $n$), and let $P'=C^{-1}PC$, then $P'$ is another m\'enage permutation and we call $P\sim P'$. In fact, $P'$ is just a relabel of $P$ where we replace $i$ with $i+a$. In the modified problem, we counted the numbers of equivalence classes from the construction above.

 For additional reading on the progression of the number of m\'enage permutations, we refer the reader to \cite{oeis}.
	
	In \S \ref{section: filling pairs and origami}, we show how m\'enage problem, especially the modified one, is associated with the problem constructing surfaces by identifying the edges of the polygon.
	\subsection{Outline}
	In \S \ref{section: filling pairs and origami}, we discuss some fundamental properties for a minimally intersecting filling coherent pair. In section 3, we give an explicit algorithm to construct such pairs on $S_g$. In \S \ref{section: filling pairs and origami}, we give an upper bound for the number using the construction from \S \ref{section: example}. In \S \ref{section: searching}, we give an alternative algorithm to search for such pairs with a program. And finally in section \S \ref{section: genus 2}, we discuss the case with $g=2$, where we will have a minimal of 4 intersections.
	\section{Minimal filling pairs and origami}
 \label{section: filling pairs and origami}
	 Suppose further $\alpha$ and $\beta$ intersect coherently, then the sides will also be coherent, and it's a $1-1$ origami with $2g-1$ squares.
	\begin{prop} \label{oppo}
		There exists a side of $\alpha$ (or $\beta$) such that it's not identified with the opposite.
	\end{prop}
	
	Now the opposite sides are in the different direction, so if every side of $\alpha$ is identified with the opposite, we have two pairs of $\alpha$: $a_1,a_1'$ and $a_2, a_2'$ adjacent to each other in the origami (note it's not necessary in the polygon), sharing the same section of $\beta$, with $b, b'$ as the image of that arc in the sides of the polygon. Then $b$ and $b'$ is in the right (or left) of $a_1$ and $a_2$. Let $a_3$ be 
	the other side of $b$ and $a_4$ be the other side of $b'$. Then since $b$ and $b'$ is actually the same in the origami, $a_3$ and $a_4$ has to be adjacent to each other. (See figure above) Since every side of $\alpha$ is identified with the opposite, we are actually playing the marking game on a circle with $n$ points on it: if the number of two vertices is adjacent, then the ``outer'' and ``inner'' vertices should also be adjacent. So it should look like an ``onion'' :$\cdots,3,1,2,4, \cdots$, but it's a circle with odd vertices, so it's impossible.
	
	\begin{figure}[h]
		\scalebox{0.4}{\includegraphics[angle=0,origin=c]{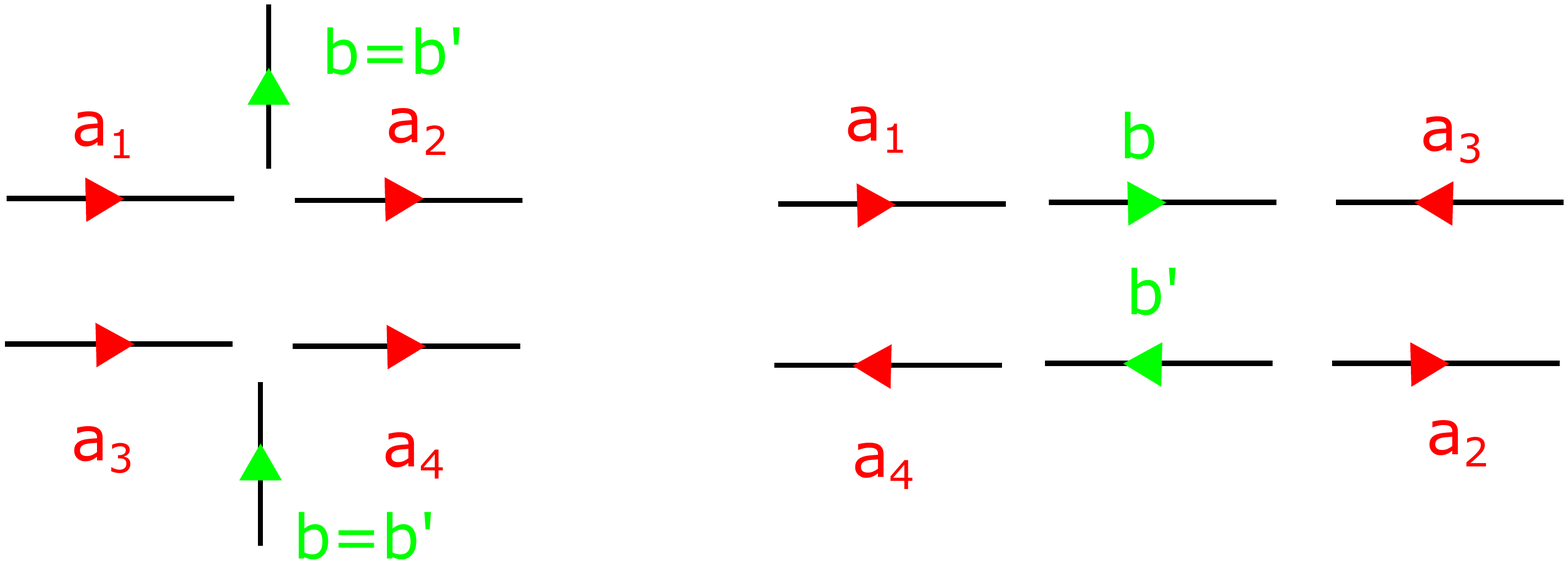}}
		\vspace{-2cm}
		\caption{One example of possible orientation, the left is part of the origami while the right is part of the polygon. Please note there are other possible orientations, but $a_3$ should always be adjacent to $a_4$ in the origami.}
		\label{orientation}
	\end{figure}
	\begin{figure}[h]
		\scalebox{0.7}{\includegraphics[angle=0,origin=c]{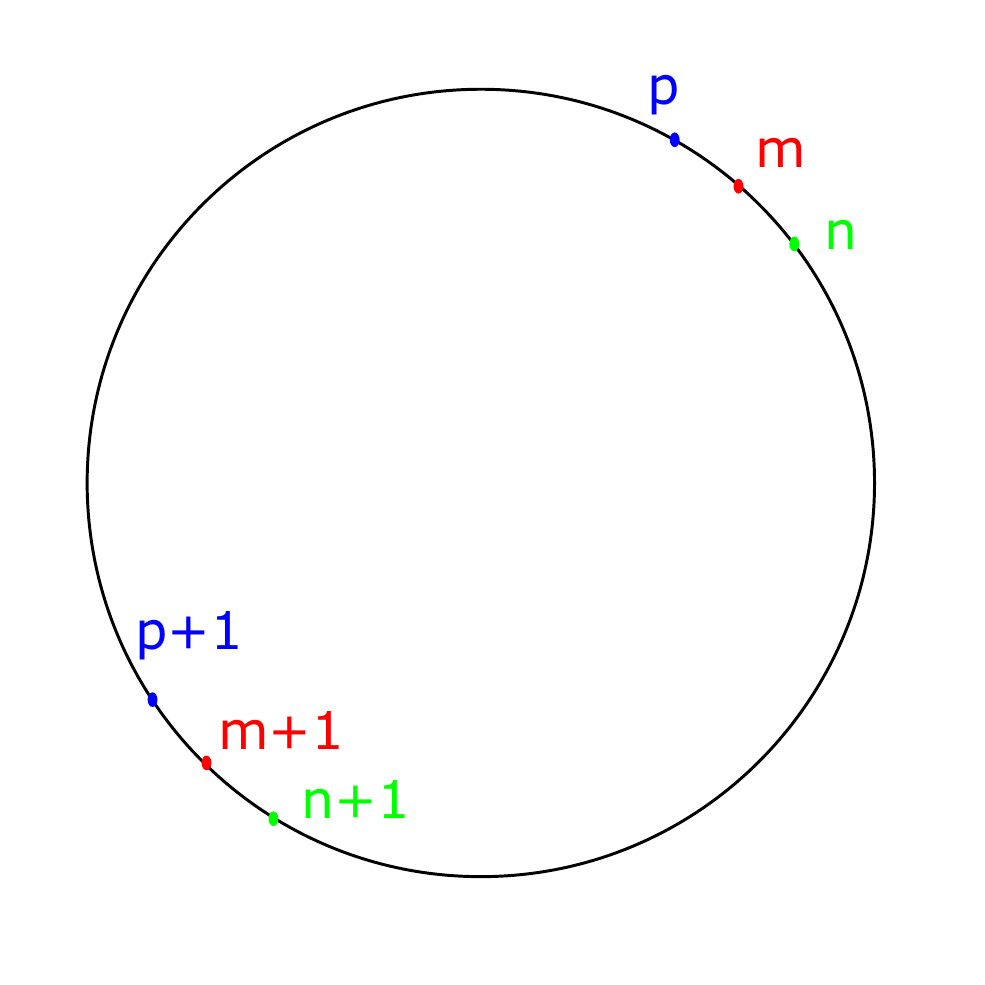}}
		\vspace{0cm}
		\caption{The game of putting numbers on vertices: if the number of two vertices (red)is adjacent, then the ``outer''  and ``inner''  vertices (blue and green) should also be adjacent. See Fig. \ref{orientation}.}
		\label{outer and inner}
	\end{figure}
	
	\begin{de}
		Let $a_i$ (i=1,2,...,2(2g-1)) to be the edges in clockwise direction from the curve $\alpha$ of the $4(2g-1)$ gon mentioned above. We define the distance between two edges, $a_i$ and $a_j$, to be$$
		d(a_i,a_j)=
		\begin{cases}
		|i-j|, |i-j|\le 2g-1\\
		2(2g-1)-|i-j|, |i-j|>2g-1
		\end{cases}$$If $a_i$ and $a_j$ ($i\not=j$) is identified in the curve $\alpha$, we define the distance of the identified segment as $d(a_i,a_j)$.
	\end{de}
	
	\begin{lemma} \label{odd}
		For the distance of two identified edges we have :\\
		(1) The distance is always odd.\\
		(2) The distance cannot be 1.
	\end{lemma}
	
	Since $\alpha$ and $\beta$ are intersecting coherently, in the polygon if the distance between two edges is even, they are in the same direction and vice versa. If the direction is same, the ``next'' edge following also have to identified with each other. With a simple induction, the distance of every pair of identified edges will be the same and even, which will be impossible with the fact that $2g-1$ is odd. If the distance between identified $a$ and $a'$ of $\alpha$ is 1, then there will be an arc $b$ of $\beta$ lies between them, which is also impossible.
	\begin{figure}[h]
		\scalebox{0.4}{\includegraphics[angle=0,origin=c]{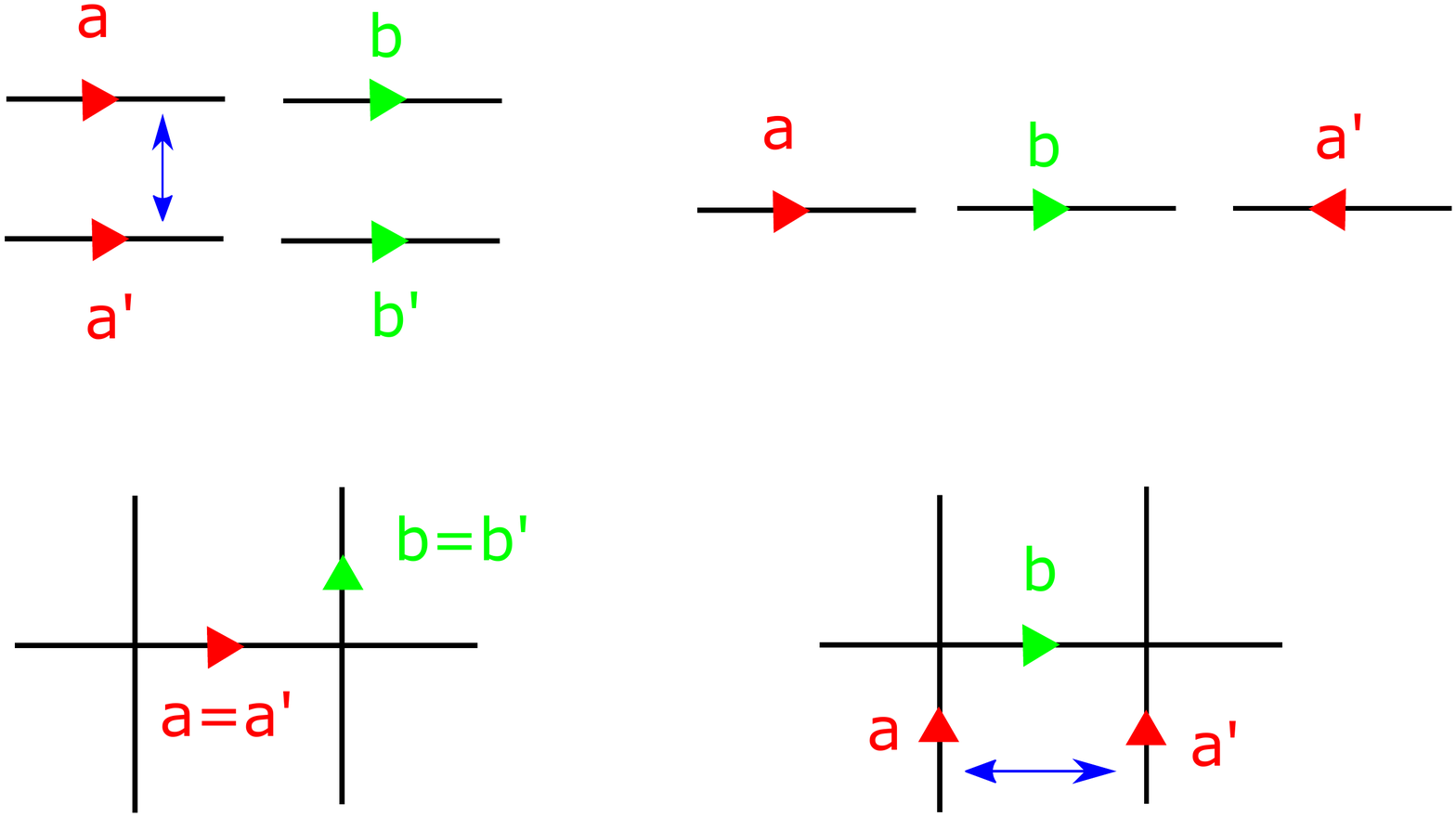}}
		\vspace{0cm}
		\caption{Left: if two edges with same direction are identified, the next one will also be identified; right: if the distance is 1, there will be one arc intersecting two identified arcs in the same direction.}
		\label{identified edges}
	\end{figure}
	
	The algorithm to transfer a polygon into an origami is simple, we are just following the way of connection for each intersection. However, things will be harder if we are transforming an origami into a polygon.

	\section{An example}
	 \label{section: example}

	We are now showing how to find minimal filling pairs for $S_g$, we take $g=3$ as an example:
	\begin{prop} \label{g3}
		There's a unique pair of filling curves in $S_3$ that is coherent with minimal intersecting number.
	\end{prop}
	
	\subsection{Step 1: Find possible ways to identify edges from curve $\alpha$ on the 20-gon.}
	
	If we check-board the edges of $\alpha$ in the $4(2g-1)$-gon (there are $2(2g-1)$ of them), according to Lemma \ref{odd}, two identified edges have to be one white and one black, and their distance cannot be $1$. And if we rotate the table and relabel them, we are actually getting the same way to identify the edges since the edges are not ordered. So, this is actually the modified m\'enage problem discussed in introduction, and as a result, when $n=5$, there are $5$ equivalent classes. However, the case that all labels are opposite is impossible due to Proposition \ref{oppo}, so there are $4$ cases in total.
	
	\begin{figure}[h]
		\scalebox{0.3}{\includegraphics[angle=0,origin=c]{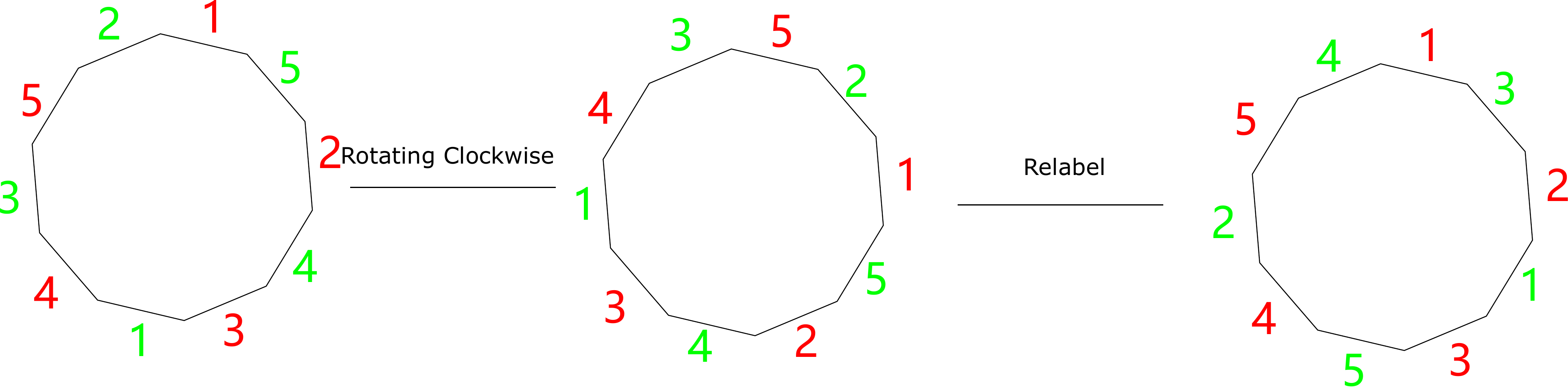}}
		\vspace{0cm}
		\caption{An example of rotating and relabelling, red are ladies that are fixed and green are gentlemen. Arcs belong to $\beta$ are not shown.}
		\label{rotation}
	\end{figure}
		
	We notice we are actually finding ``permutations'' of the m\'enage problem, which is introduced in \cite{gen}. According to the Table.1 of that paper, and notice Proposition \ref{oppo}, the case that all couples are seated opposite is impossible. So we have:
	
	\begin{prop}
		The only possible ways to identify edges of $\alpha$ will be one of the four following cases up to symmetry.
	\end{prop}
	\begin{figure}[H]
		\scalebox{0.3}{\includegraphics[angle=0,origin=c]{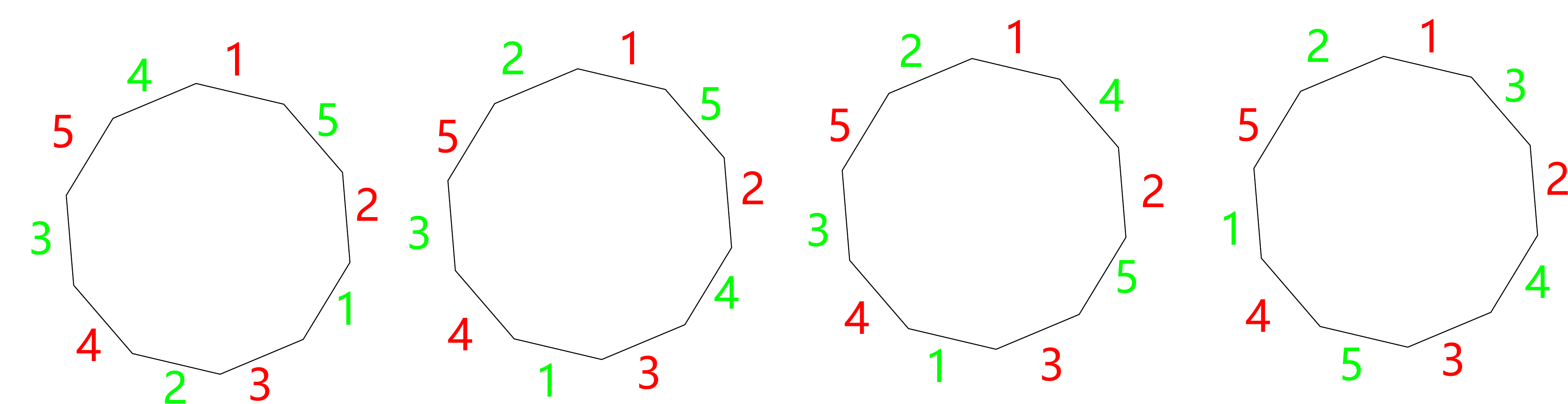}}
		\vspace{0cm}
		\caption{Four possible cases. Arcs belong to $\beta$ are not shown.}
		\label{4 possible cases}
	\end{figure}
	\subsection{Step 2: Find possible labels of edges in $\beta$ for edges in $\alpha$}

\begin{figure}[H]
	\scalebox{0.4}{\includegraphics[angle=0,origin=c]{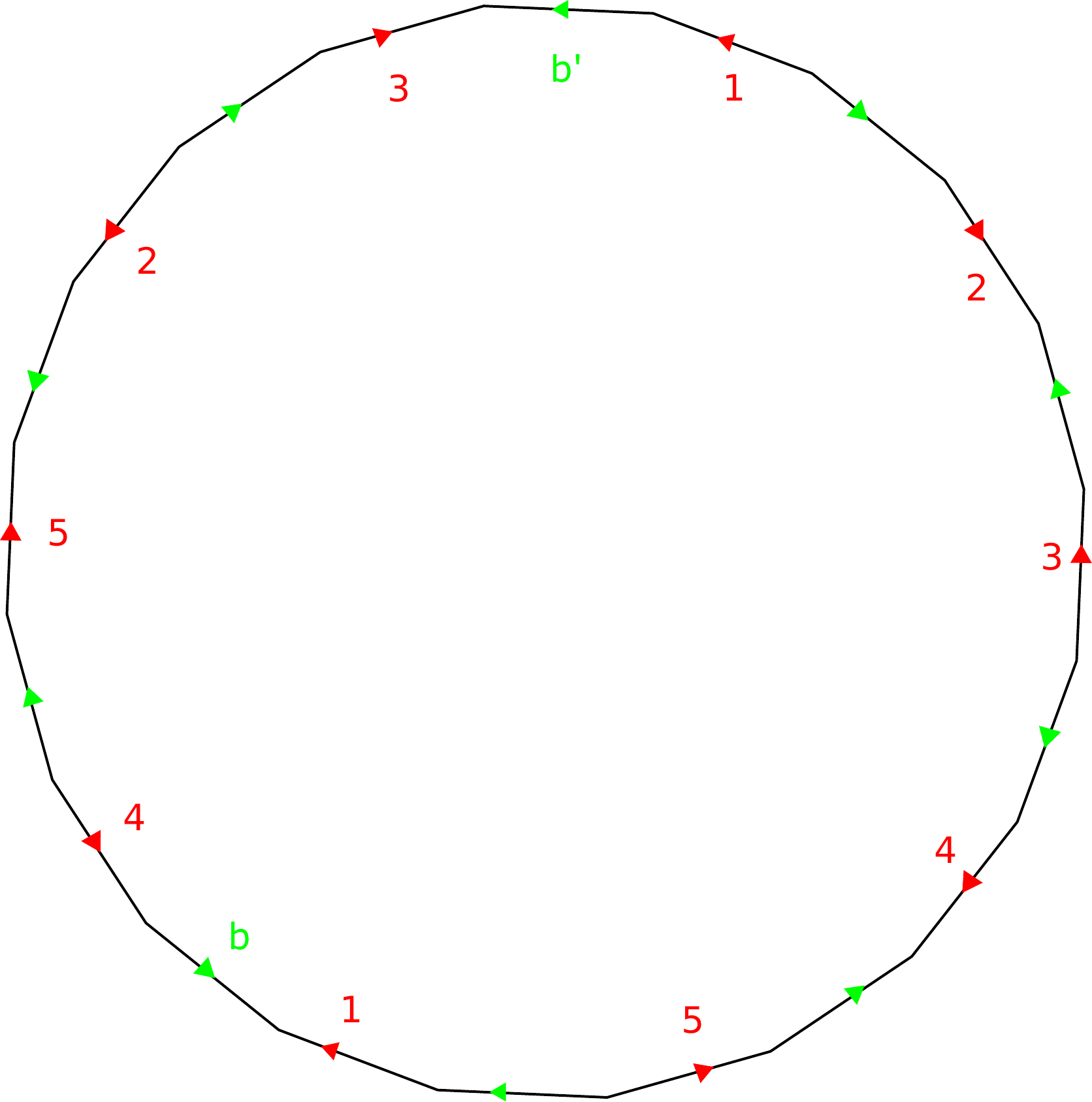}}
	\vspace{0cm}
	\caption{Putting the labels of $\alpha$ on the polygon.}
	\label{loop: step 1}
\end{figure}

	We give an orientation of $\beta$ and label the arcs in $\beta$ with the orientation. Suppose the curves are orientated as the left of the picture, and two $a$s are identified as an arc of $\alpha$ like the middle, then the right shows it will behave like that in the origami and thus the label of $b'$ has to be 1 bigger than $b$.
	We take the loop (1 2 3 4 5 1 4 5 2 3) as an example (the reader can check this is a just a relabel of the second case in \ref{4 possible cases}), first we put the labels of $\alpha$ on the polygon, see Figure \ref{loop: step 1}
	
	According to the discussion above, we notice that the label of $b'$ have to be the next of $b$. Without loss of generality, we let the label of $b$ to be 1, so we continue this process until we label five edges according to one direction of $\alpha$:
	
\begin{figure}[H]
	\scalebox{0.4}{\includegraphics[angle=0,origin=c]{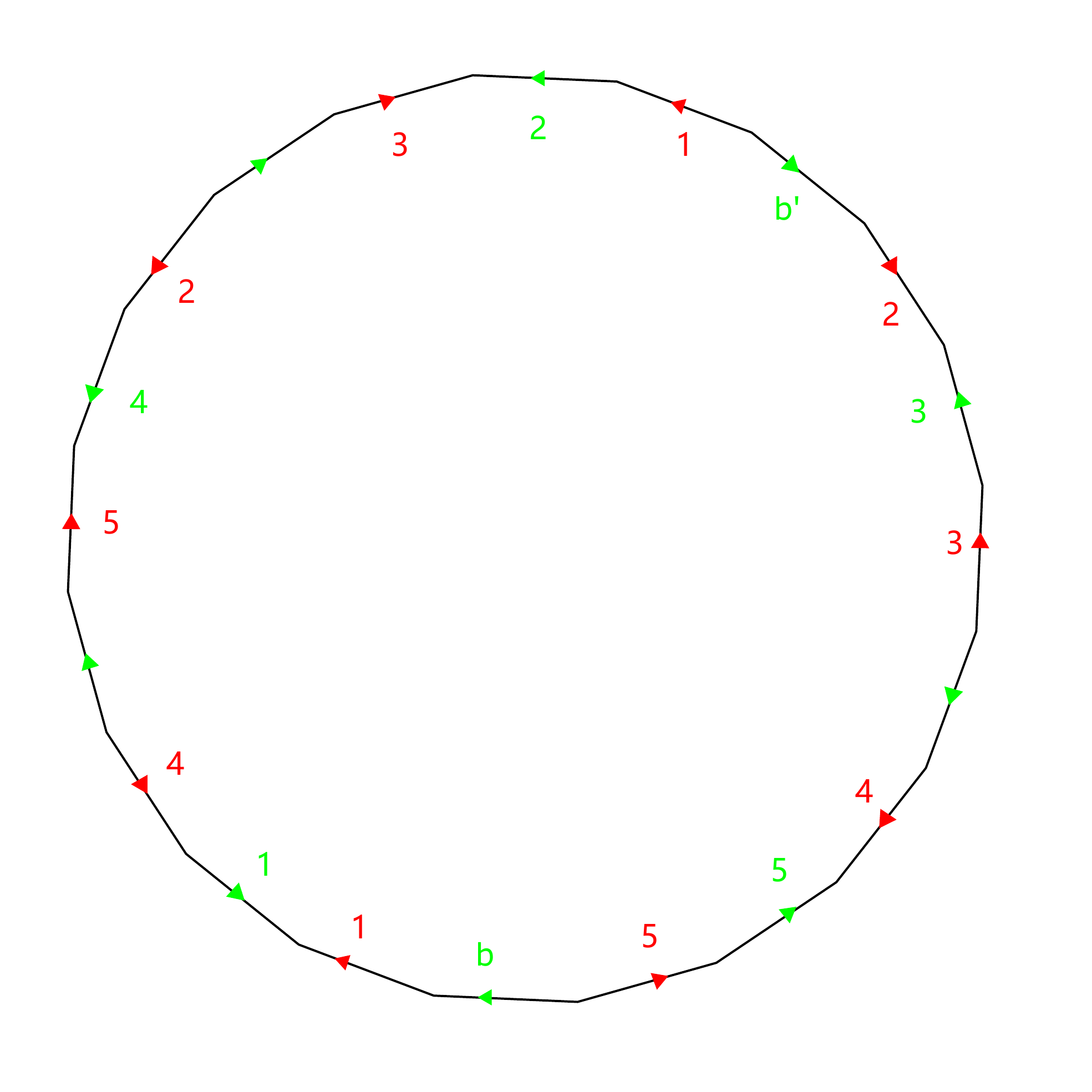}}
	\vspace{0cm}
	\caption{Labelling five edges according to one direction of $\alpha$.}
	\label{loop: step 2}
\end{figure}
	
	We continue on labelling the other sides, notice $b'$ have to be the next of $b$, however since the other direction is already labelled, we label $b$ as $x$, so $b'$ will be $x+1$.
	
\begin{figure}[H]
	\scalebox{0.4}{\includegraphics[angle=0,origin=c]{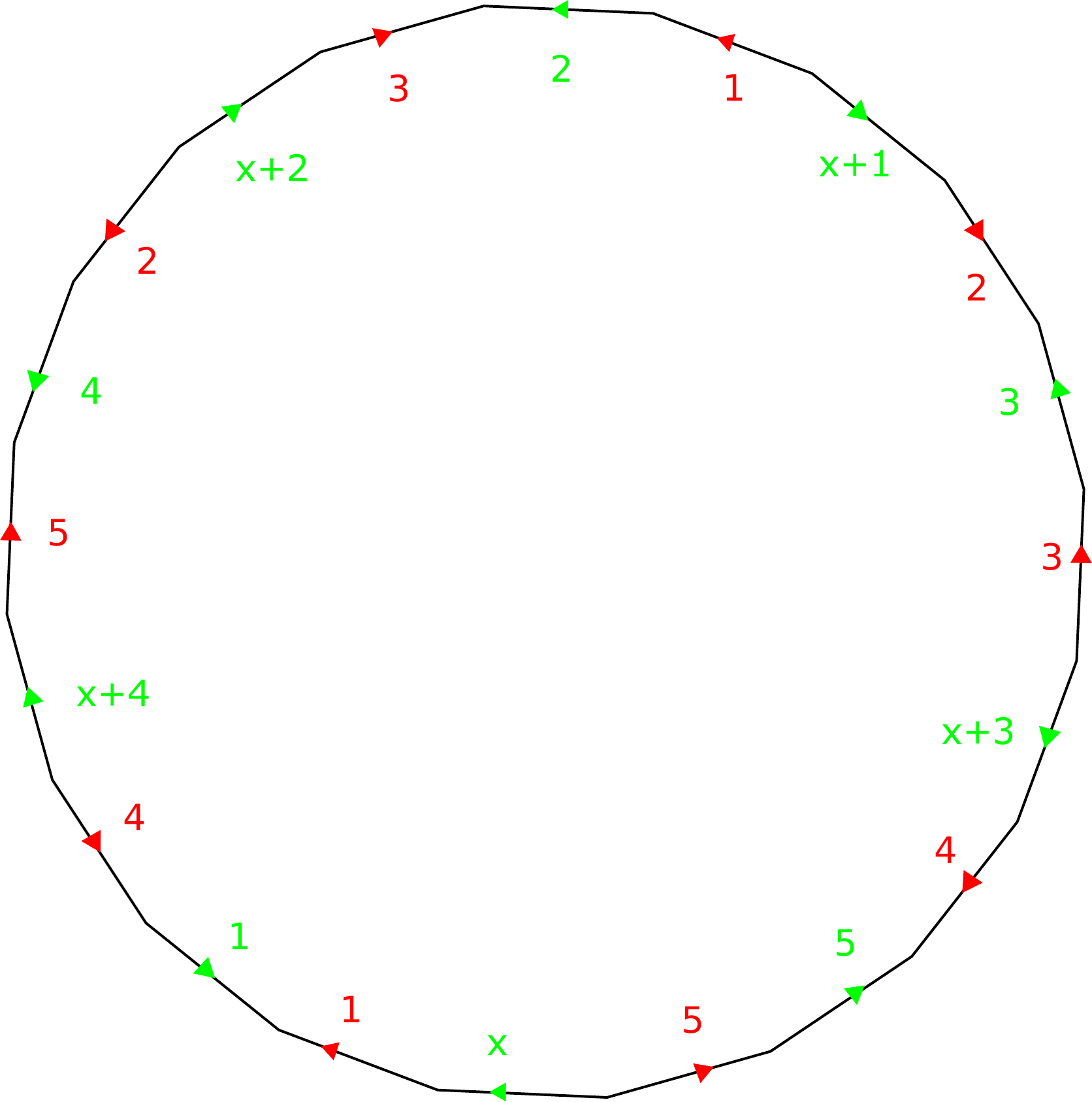}}
	\vspace{0cm}
	\caption{Labelling the other sides.}
	\label{loop: step 3}
\end{figure}

	Similarly, we finished labelling all the edges in $\beta$.
	
	\subsection{Step 3: Identify edges in $\beta$.}
	
	We want to find what $x$ is in the last step, where $x$ can be an integer from 1 to 5. Notice that in the lower half of the origami and the upper half origami will be like the following:
	
\begin{figure}[H]
	\scalebox{0.5}{\includegraphics[angle=0,origin=c]{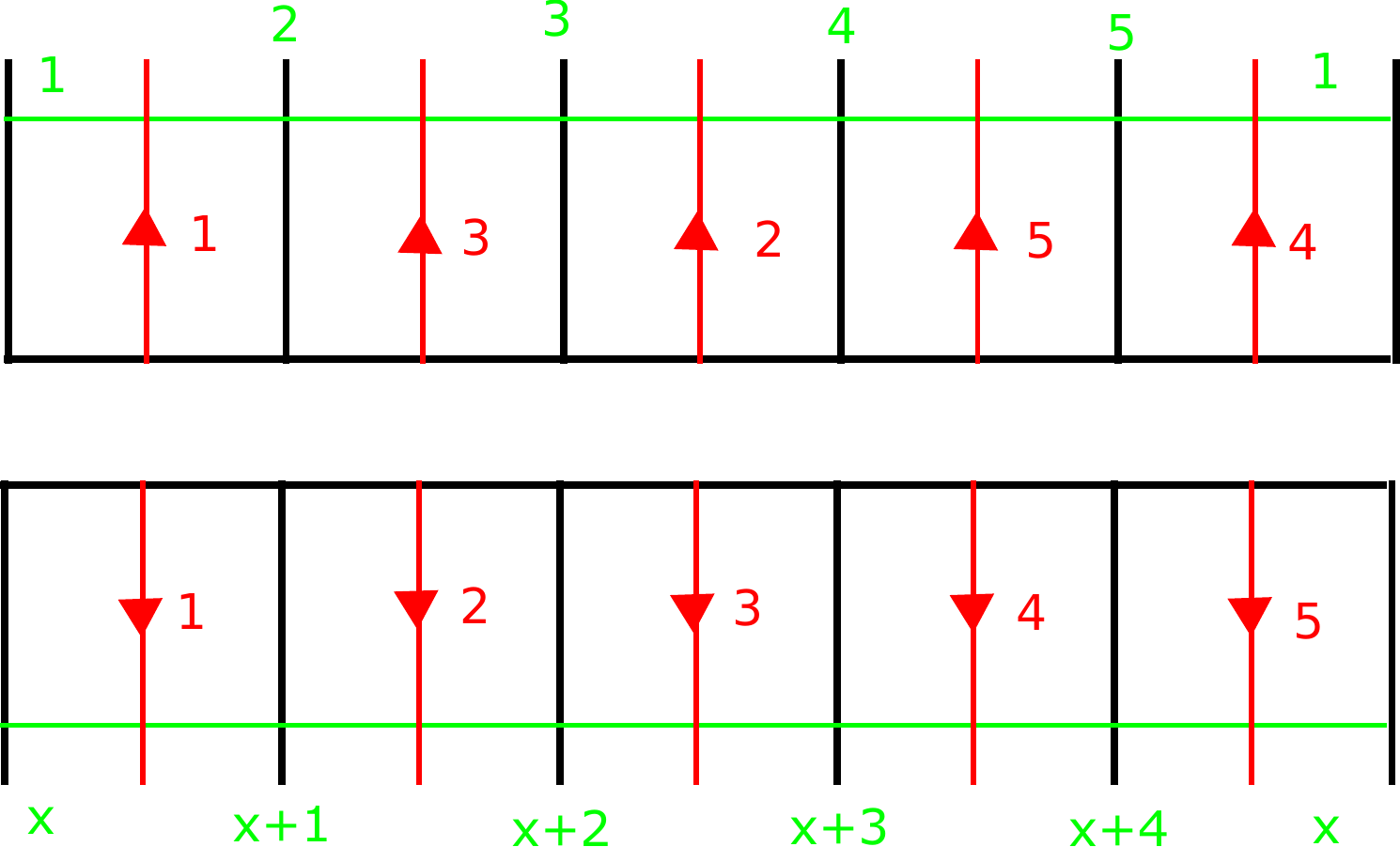}}
	\vspace{0cm}
	\caption{Identifying edges in $\beta$.}
	\label{identified edges in example}
\end{figure}
	
	We say $p$ is \emph{the permutation of the $1-1$ origami} if the top edge of $i$-th square is identified with the bottom edge $p(i)$-th square. 
	
	So if $x=1$, it is not an origami since edge 1 in $\alpha$ is followed by itself, same happen with $x=2$ where edge 2 is followed by itself, and when $x=5$, edge 3 is followed by itself. When $x=3$, it's an origami with permutation (1 2 5 3 4), and when $x=4$, it's an origami with permutation (1 5 2 4 3). For (1 2 5 3 4), the different between two adjacent entries are $1,3,3,1,2$ while for (1 5 2 4 3) they are $-1,-3,-3,-1,-2$ so they are just the mirror reflection of each other, so there's one (up to relabelling) possible labelling of the 20-gon with the loop (1 2 3 4 5 1 4 5 2 3). Similarly we can prove there's no possible labelling with the other three loops and this finish the proof of the proposition. 
	
	We can use similar algorithm on higher genus surfaces, however the way to identify $\beta$ may not be unique. For example, in $S_5$, if the top of $\beta$ is $1,2,9,3,4,6,7,5,8$ and the bottom is $x,x+1,x+9,x+3,x+4,x+6,x+7,x+5,x+8$, then when $x=8$ and when $x=1$, they will form two origamis with permutation (1 2 4 5 8 6 3 7 9) and (1 3 9 7 6 8 5 4 2). Taking the different between two adjacent entries we get $1,2,1,3,7,1,2,3,7$ and $-7,-3,-2,-1,-7,-3,-1,-2,-1$. They are neither the same nor mirror reflection.
	
	\section{An estimate on upper bound for the \\ number of ordered filling pairs}
	 \label{section: estimate}
	
	Notice in Step 1, according to \cite{oeis}, the number of m\'enage permutations $h(g)\sim\frac{(2g-2)!}{e^2}$.
	
	In Step 2, for fixed $\alpha$, the $\beta$ in the two half squares are uniquely determined.
	
	In Step 3, when all edges in $\beta$ are labelled, we may have up to $2g-2$ ways to identify them and they are actually symmetric (notice 1 cannot be identified with 1).
	
	So according to the steps, notice we counted twice for mirror reflections, there will be at most $(2g-2)h(g)/2$ possibilities and this finish the proof of Theorem \ref{main}.
	\section{Searching for possible origamis}
	 \label{section: searching}
	
	In the discussion above, given an origami with $2g-1$ squares, we have an algorithm to find the arcs of $\alpha$ on the corresponding $[4(2g-1)]$-gon. Suppose the origami is already oriented, we start from an arc of $\alpha$ in the origami and locate it on the polygon. According to the orientation of the polygon and the origami, we can find the location on the polygon for the adjacent arc of $\alpha$. Then we repeat this process until all the polygon is filled. Notice the orientation of the polygon is coherent, so if the last step is finding the right neighbour of an arc on the upper half of the origami, the next will be finding the left neighbour of the new arc on the lower half, and vice versa.
	
	Take the origami below as example, we start from arc 1 in $\alpha$ and locate it in one of the arcs in polygon as shown in the picture. We notice the two blue areas in the origami and the polygon have the same direction, so we move one grid right in the lower half of the origami and find the next arc is 4. Again, we find that the direction for two light blue area in the origami and the polygon are also the same, so we move one grid left in the upper half of the origami and the result is 5.
	
	\begin{figure}[h]
		\scalebox{0.4}{\includegraphics[angle=0,origin=c]{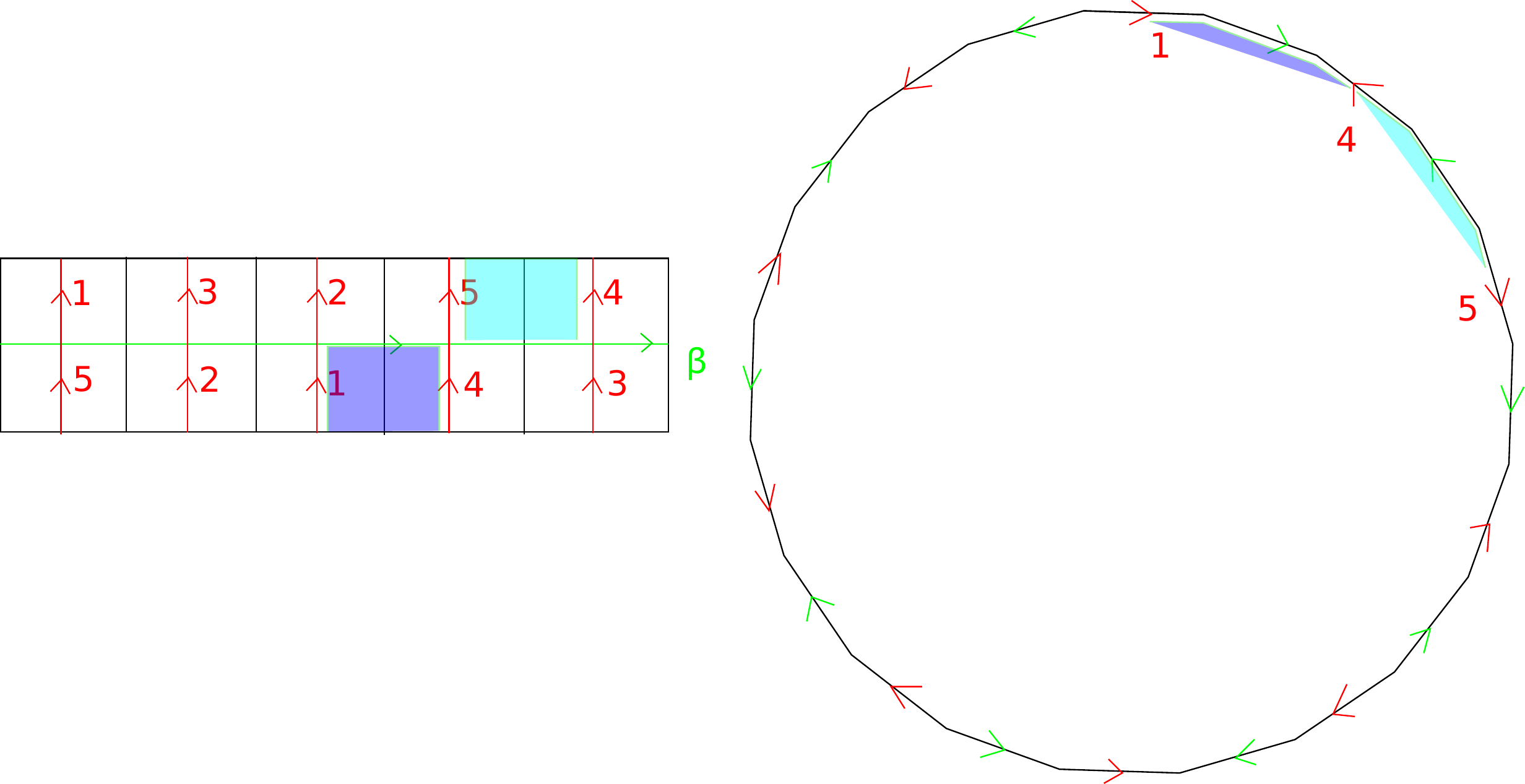}}
		\vspace{0cm}
		\caption{Example of an origami and its corresponding polygon.}
		\label{search}
	\end{figure}
	
	\begin{prop}
		An origami determines a minimal filling pair if and only if they can fill the arcs of $\alpha$ in the $4(2g-1)-gon$ with the algorithm above and the distance of two identified edges are always odd.
	\end{prop}
	
	With the proposition above, we have an algorithm to check whether an origami determines a minimal filling pair. And if we check all the possible $1-1$ origamis with $2g-1$ squares, we will be able to find all ordered minimal filling pairs. Notice if we change the start of $\beta$, we will have $2g-1$ different origamis but they all determine the same origami; and if we change the orientation of $\beta$, we will find 2 origamis. So each ordered minimal filling pair will be correspond to $2(2g-1)$ origamis. With the discussion above I made a searching program \url{https://github.com/expectedid/countorigami}. Please notice Figure 3 in \cite{AMN} counted the number of $SL(2, \mathbb{Z})$ orbits which include shearing and rotating but my program is counting the number of all possible origamis up to mirror symmetry and relabelling (no shearing), so the result may be different and mine will be larger for the same genus.
	
	The results with smaller genera are:
 \begin{center}
	\begin{tabular}{cccc}
		\hline
		\makecell[c]{Genus of\\the surface}&\makecell[c]{Number of squares\\in origami}&\makecell[c]{Number of\\coherent minimal\\filling pairs}&$[\frac{(g-1)(2g-2)!}{e^2}]$\\
		\hline
		3&5&1&6\\
		4&7&8&292\\
		5&9&436&21826\\
		6&11&23904&2455523\\
		7&13&2448720&388954903\\
		\hline
	\end{tabular}
 \end{center}
	
	\section{Case with genus 2}
	 \label{section: genus 2}
	 
	When $g=2$, \cite{AH} suggests that for the number of squares, the bound $2g-1=3$ cannot be obtained and instead the minimal number of squares for the filling pairs is 4 with 2 boundary components in $S_2-(\alpha\cup\beta)$. There are, $(4-1)!=6$, $1-1$ origamis with four squares and we check the Euler characteristic for all of them. As a result, there are 4 out of 6 whose Euler characteristic is $-2$, which means the surface is $S_2$. However, up to symmetry, they are actually the same pair of curves on the surface. So, we can improve the proposition into following:

	\begin{thm}
		There's a unique pair of minimal coherent filling curves in $S_2$ or $S_3$ that is coherent with minimal intersecting number.
	\end{thm}
	
	We have discussed the genus $3$ case in \ref{g3}. Figure below is one of the genus $2$ origami and all genus $2$ origami is in the same $Mod^{\pm}(S_2)$ orbit.
	
		\begin{figure}[h] 
		\scalebox{1}{\includegraphics[angle=0,origin=c]{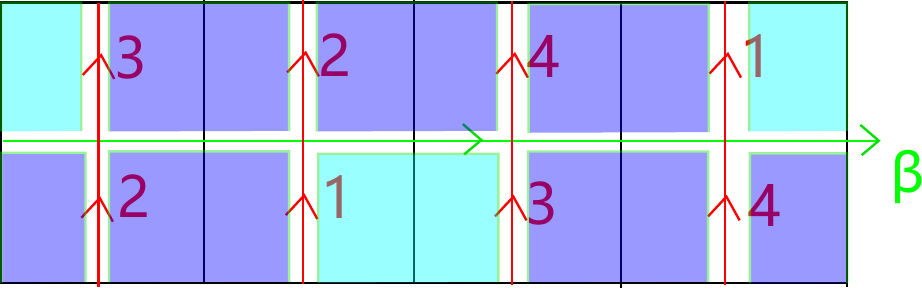}}
		\vspace{0cm}
		\caption{Example of an genus 2 origami, the violet is one disk in $S_2-(\alpha\cup\beta)$ while the blue one is the other.}
	\end{figure} \label{g2}
	
\end{document}